\documentclass[preprint,12pt, a4paper]{elsarticle}

\usepackage{amssymb}
\usepackage{hyperref}

\usepackage{listings}
\usepackage{color, colortbl}

\newcommand{\athena}{\textsc{ATHENA} }
\newcommand{\TODO}[1]{{\color{red}#1}}

\journal{Software Impacts}

\begin{document}

\begin{frontmatter}



\title{ATHENA: Advanced Techniques for High Dimensional Parameter Spaces to
  Enhance Numerical Analysis}


\author[]{Francesco Romor}\ead{francesco.romor@sissa.it} \author[]{Marco
Tezzele}\ead{marco.tezzele@sissa.it} \author[]{Gianluigi
Rozza}\ead{gianluigi.rozza@sissa.it}

\address{Mathematics Area, mathLab, SISSA, International School of Advanced Studies, Trieste, Italy}

\begin{abstract}
  \athena is an open source Python package for reduction in parameter space. It
  implements several advanced numerical analysis techniques such as Active
  Subspaces (AS), Kernel-based Active Subspaces (KAS), and Nonlinear Level-set
  Learning (NLL) method. It is intended as a tool for regression, sensitivity
  analysis, and in general to enhance existing numerical simulations' pipelines
  tackling the curse of dimensionality.
\end{abstract}

\begin{keyword}
Parameter space reduction \sep Active subspaces \sep Kernel-based active
  subspaces \sep Nonlinear level-set learning \sep Python



\end{keyword}

\end{frontmatter}

\section{The \athena package}
Parameter space reduction techniques are particularly suited to tackle the curse
of dimensionality affecting high dimensional problems for many-queries and
outer-loop applications such as optimization, uncertainty propagation, and
statistical inference. The increased capabilities in simulating every day more
and more complex phenomena with many input design parameters need to be
compensated by a proper handling of the dimension of such parameters. This is
due to the fact that the amount of samples needed scales superlinearly with
respect to the input dimension. Possible strategies to mitigate this issue could
be sparse grid sampling, adaptive sampling, and model order
reduction~\cite{morhandbook2019, rozza2018advances}.

\begin{figure}[htb]
  \centering
  \includegraphics[width=.85\textwidth]{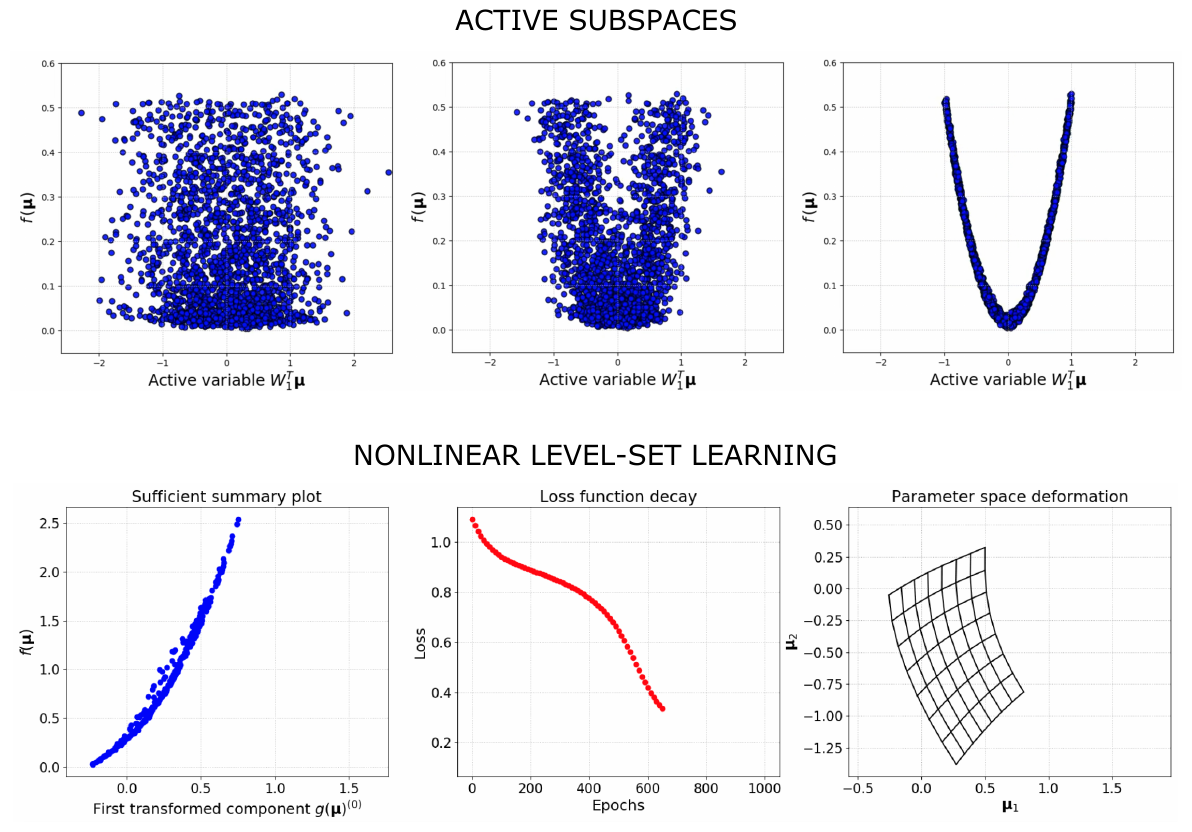}
  \caption{Illustration of the application of active subspaces (top panel) and
    nonlinear level-set learning (bottom panel).}
  \label{fig:athena_scheme}
\end{figure}

In this package we focus on techniques which use linear and nonlinear
transformations to align the input space along the directions of maximum
variation of the function of interest. In particular in the
 \athena open source Python
package are implemented the interfaces to easily employ Active
Subspaces~\cite{constantine2015active}, Kernel-based Active
Subspaces~\cite{romor2020kas}, and the Nonlinear Level-set Learning
\cite{zhang2019learning} method. Figure~\ref{fig:athena_scheme} depicts a
simple application of AS and NLL.

\begin{figure}[htb]
  \centering
  \includegraphics[width=.85\textwidth]{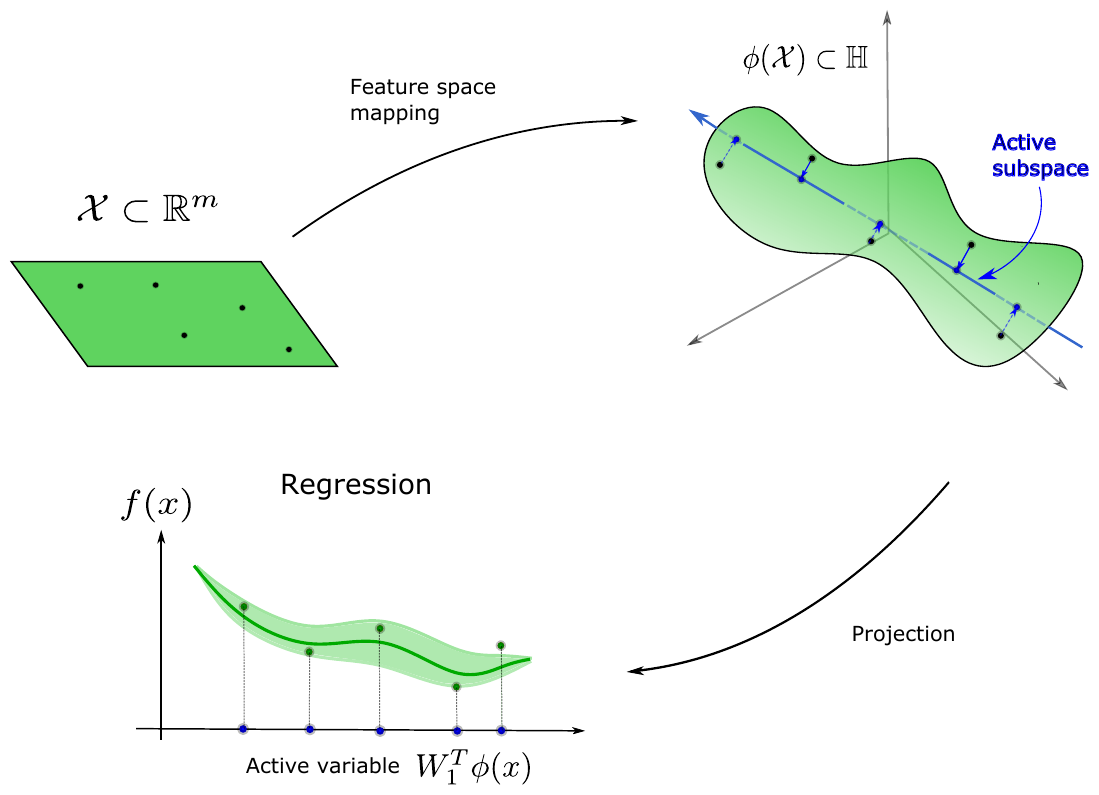}
  \caption{Illustration of the application of kernel-based active subspaces.}
  \label{fig:kas_scheme}
\end{figure}

These techniques are purely data driven, since they necessitate only
input-output couples. We emphasize that even if the methods are gradient-based,
it is possible to approximate the gradients of the target function with respect
to the parameters using only input-output data. Thus they are particularly
suited to enhance existing numerical simulations' pipelines such as optimization
procedures, possibly including commercial softwares. Moreover it is possible to
couple parameter space reduction with reduced order methods, enabling higher
speedup computations with controlled accuracy, especially for the solution of
parametric PDEs, both in academia and in industry.

The package is multi-purpose and intended for practitioners, engineers, data
scientists, and researchers.

\section{The impact to research fields}
The package made possibile the creation of several new methods for
optimization, regression, and classification tasks.
The interface for Active Subspaces studies of \athena is used
in~\cite{demo2020asga} to enhance genetic algorithms when optimizing complex and
high-dimensional functions. The reduction affects the dimension of the
population of the genetic algorithm at each step. The same interface for Active
Subspaces is used in~\cite{romor2020pamm} to improve the approximation of
Gaussian process regression of scalar functions with low intrinsic
dimensionality. The information regarding the presence of an Active Subspace is
delivered to a nonlinear multi-fidelity model thus able to reach a higher
accuracy. In some cases the Active Subspace is not present due to its linear
nature, an extension is introduced in~\cite{romor2020kas} with Kernel-based
Active Subspaces. In this case the inputs are projected to a high-dimensional
space and then reduced, mimicking the framework of reduction in Reproducing
Kernel Hilbert Spaces. Figure~\ref{fig:kas_scheme} shows a sketch to
understand how the KAS method works. In the context of a computational
fluid dynamic benchmark, \textsc{ATHENA}'s Kernel-based Active Subspaces interface is tested on the
design of a response surface showing an improvement of the accuracy with respect
to Active Subspaces.
In~\cite{romor2021las} a local approach to parameter space reduction
is presented to improve the accuracy of regression and classification.

\section{The impact to industrial collaborations}
\athena has been created thanks to an industrial collaboration with
Fincantieri S.p.A., matching the demand of an intuitive, easily
integrable tool for parameter space reduction in the framework of
structural optimization of modern passenger ship hulls. In an ongoing
industrial Ph.D. grant sponsored by Electrolux Professional, the
package is used to reduce the dimensionality of convolutional
artificial neural networks for image recognition.

\athena is also utilized in~\cite{tezzele2020enhancing} to reduce the dimensionality
of a response surface resulting from the shape optimization of an airfoil. The
proposed pipeline couples Dynamic Mode Decomposition with Dynamic Active
Subspaces to reduce the overall computational resources. The methods are used
also in~\cite{demo2021hull} to optimize the shape of a benchmark hull using
advanced geometrical mophing.

\section{Declaration of competing interest}
The authors declare that thay have no known competing finantial interests or
personal relationships that could have appeared to influence the work reported
in this paper.

\section*{Acknowledgements}
This work was partially supported by an industrial Ph.D. grant sponsored by
Fincantieri S.p.A. (IRONTH Project), by MIUR (Italian ministry for university
and research) through FARE-X-AROMA-CFD project, and partially funded by European
Union Funding for Research and Innovation --- Horizon 2020 Program --- in the
framework of European Research Council Executive Agency: H2020 ERC CoG 2015
AROMA-CFD project 681447 ``Advanced Reduced Order Methods with Applications in
Computational Fluid Dynamics'' P.I. Professor Gianluigi Rozza.





\bibliographystyle{elsarticle-num}
\bibliography{athena_bib.bib}

\section*{Illustrative Examples}
Optional : you may include one explanatory  video that will appear next to your
article, in the right hand side panel. (Please upload any video as a single
supplementary file with your article. Only one MP4 formatted, with 50MB maximum
size, video is possible per article. Recommended video dimensions are 640 x 480
at a maximum of 30 frames / second. Prior to submission please test and validate
your .mp4 file at
\url{http://elsevier-apps.sciverse.com/GadgetVideoPodcastPlayerWeb/verification}
. This tool will display your video exactly in the same way as it will appear on
ScienceDirect. )

\section*{Required Metadata}
\label{}

\section*{Current code version}
\label{}

Ancillary data table required for subversion of the codebase. Kindly replace
examples in right column with the correct information about your current code,
and leave the left column as it is.

\begin{table}[!h]
\begin{tabular}{|l|p{6.5cm}|p{6.5cm}|}
\hline
\textbf{Nr.} & \textbf{Code metadata description} & \textbf{Please fill in this
column} \\
\hline
C1 & Current code version & v0.1.2 \\
\hline
C2 & Permanent link to code/repository used for this code version & \url{https://github.com/mathLab/ATHENA} \\
\hline
C3  & Permanent link to Reproducible Capsule & \TODO{Add repr capsule}\\
\hline
C4 & Legal Code License   & MIT License (MIT) \\
\hline
C5 & Code versioning system used & git \\
\hline
C6 & Software code languages, tools, and services used & Python \\
\hline
C7 & Compilation requirements, operating environments \& dependencies
                                                  & ATHENA requires
                                                    \texttt{numpy},
                                                    \texttt{matplotlib},
                                                    \texttt{scipy},
                                                    \texttt{torch},
                                                    \texttt{GPy},
                                                    \texttt{GPyOpt},
                                                    \texttt{scikit-learn},
                                                    \texttt{sphinx}
                                                    (for the
                                                    documentation) and
                                                    \texttt{nose}(for
                                                    local test)\\
\hline
C8 & If available Link to developer documentation/manual & \url{https://mathlab.github.io/ATHENA/} \\
\hline
C9 & Support email for questions & marco.tezzele@sissa.it francesco.romor@sissa.it\\
\hline
\end{tabular}
\caption{Code metadata (mandatory)}
\label{}
\end{table}

\section*{Current executable software version}
\label{}

Ancillary data table required for sub version of the executable software: (x.1,
x.2 etc.) kindly replace examples in right column with the correct information
about your executables, and leave the left column as it is.

\begin{table}[!h]
\begin{tabular}{|l|p{6.5cm}|p{6.5cm}|}
\hline
\textbf{Nr.} & \textbf{(Executable) software metadata description} &
\textbf{Please fill in this column} \\
\hline
S1 & Current software version & v0.1.2 \\
\hline
S2 & Permanent link to executables of this version  & \url{https://github.com/mathLab/ATHENA/releases/tag/v0.1.2} \\
\hline
S3  & Permanent link to Reproducible Capsule & \TODO{Add capsule}\\
\hline
S4 & Legal Software License & MIT License (MIT) \\
\hline
S5 & Computing platforms/Operating Systems & Linux, OS X, Unix-like
\\
\hline
S6 & Installation requirements \& dependencies & ATHENA requires
                                                 \texttt{numpy},
                                                 \texttt{matplotlib},
                                                 \texttt{scipy},
                                                 \texttt{torch},
                                                 \texttt{GPy},
                                                 \texttt{GPyOpt},
                                                 \texttt{scikit-learn},
                                                 \texttt{sphinx} (for
                                                 the documentation)
                                                 and \texttt{nose}(for
                                                 local test) \\
\hline
S7 & If available, link to user manual - if formally published include a
reference to the publication in the reference list & \url{https://mathlab.github.io/ATHENA/} \\
\hline
S8 & Support email for questions & marco.tezzele@sissa.it francesco.romor@sissa.it\\
\hline
\end{tabular}
\caption{Software metadata (optional)}
\label{}
\end{table}

\end{document}